\documentclass[12pt]{article}
\title{On Two Problems From "Hyperidentities and Clones"} 

\author{Gerhard R. Paseman\\
International Congress of Mathematicians 2014, Logic Poster Section 01} 
\date{Seoul, Republic of Korea\\
August 14, 2014}
\usepackage[left=1cm,top=1cm,right=1cm,bottom=1cm,bindingoffset=0cm]{geometry}
\usepackage[utf8]{inputenc}
\usepackage{CJKutf8}
\usepackage[some]{background}
 
\newcommand\textkorean[1]{%
  \begin{CJK}{UTF8}{mj}#1\end{CJK}}

\begin{document}
\newcommand{\ideq}{\approx}
\newcommand{\myfill}{\vspace{1cm}}

See second panel (page 2?) for title and abstract.

\section*{Poster construction}

There are twelve panels, arranged in four rows
of three panels in each row to be placed on
a background with panel ratio 4:3.  Final size will
be close to the paper size + 5\% in both dimensions,
but limiting space is 1.6 meters tall and .9 meters
wide.  This page has meta information and is not part of
the eventual display.

The assembly will be folded and unfolded, so each
panel will have some holes cut out of the corners
to ease folding.  There will be an alternating
pattern of circle and diamond shapes, with a circle
at the top left between the upper left and upper
middle panel, and a diamond at the top right
between the upper middle and upper right panel.
The holes will have a radius about 5\% of the height
of the panel.
(Yes, you heard right: the presentation has holes
in it.)

The paper will be centered on each panel, and
likely have the same aspect ratio.  At present the
paper size is 14 inches by 10.5 inches.

In addition to the text, there will be a watermark
pattern underlying the text. (At present, I cannot use
the ArXiv system to display the characters.  I apologize
for the inconvenience.)  
There will be a Korean word for logic (nonli \textkorean{??}),
one for mathematics (suhag \textkorean{??}) and one for
welcome (hwan-yeong \textkorean{??}), 
 as well as ICM and 2014 in some
arrangement.  The current arrangement is to have
ICM vertically on the right most upper three panels,
2014 on the bottom three panels, and the remaining
panels have welcome on the top row, logic on the
second row, and mathematics on the third row.
(At this time, the watermark/panel arrangement is
\begin{tabular}{ccc}
hwan- & yeong & I \\
non   & li    & C \\
su    & hag   & M \\
2 & 01 & 4 \\
\end{tabular}

I thank Joy Song for her consulting time on the
arrangement and choice of Korean words to use.

Note: to fix margins, I have placed what is
geometrically the first (hwan-) panel second in this
sequence.  In assembling, be sure to arrange the
hwan- and yeong panels so that the title (yeong panel)
 is in the middle.

\newpage
\SetBgColor{red}
\SetBgOpacity{0.2}
\SetBgAngle{0}
\SetBgScale{50}
\SetBgHshift{-6.25}
\SetBgContents{\textkorean{??}}
\BgThispage

\maketitle

\abstract{A hyperidentity $E$ can be viewed as a statement in second order logic.
When combined with a similarity type $\tau$, it can also be considered as
a set of first order statements. Based on examples from \cite{P1}, which included
hyperassociativity and $\tau=<2>$, it was conjectured that
each first order theory so produced was finitely axiomatizable. Part
of the analysis suggested further investigating the relatively free
2-generated semigroup satisfying one or both of the equations
$xxyxxyz=xxyyz$ and $zyyxx=zyxxyxx$.

At ICM 1994, the conjecture above was refuted, and a finite basis
problem arose: Is it decidable which pairs $<E,\tau>$ give rise to
finitely axiomatizable theories? This problem will be examined, and
its connections to other fields (e.g. symbolic dynamics) will be
reviewed. In doing so, we give partial solutions to problems
27 and 28 from \cite{DW}.}

\myfill

\newcommand{\ELL}{L}
In 1993, we came across the question:
"How small can an equational basis for hyperassociativity be?"
"Hyperassociativity" then meant
an infinite set of equations, one for each term $t$ of arity two
formed in the language $\ELL$ of one binary operation symbol, that said
that $t$ obeyed the associative law. If $\Omega_2$ is the collection
of these terms, one can write informally

$$ F(x,F(y,z)) \ideq F(F(x,y),z) \textrm { means }
\forall t \in \Omega_2 \forall xyz [t(x,t(y,z)) = t(t(x,y),z)]$$
or the \textit{hyperidentity} given by $ F(x,F(y,z)) \ideq F(F(x,y),z) $
is \textit{represented} or realized (in $\ELL$) by the infinite set
of identities. Any algebra with a binary operation that satisfied this
infinite set would be hyperassociative: its operation and all
binary derived operations would be associative.  Thus,
such an algebra would be a semigroup that obeyed the semigroup laws
$x^2=x^4$ (corresponding to $t(x,y) =x^2$), $xyxzxyx=xyzyx$,
$x^2y^2z=x^2yx^2yz$, and $xy^2z^2=xyz^2yz^2$ (corresponding to
$t(x,y)$ being $xyx$ or $xxy$ or $xyy$ respectively).

Denecke and Koppitz had investigated the question in \cite{DK} and found
a finite basis of roughly 1000 equations for this infinite set.  Independently
(and with different methods) Polak \cite{P8} and we \cite{P1} showed 
that there was a basis of five equations: $(xy)z = x(yz)$ and the four
parentheses-saving semigroup laws above.  This led to some
further investigations in structure \cite{P9} and language \cite{P2} \cite{P3} .

\myfill

These studies led to the Main Question below, which we call
the \textbf{Finite Basis Problem}: informally,
 is there a computer program that will tell us if there is a finitely based
equational theory equivalent to the one represented by the input pair $<E,\tau>$?

A more general version of this question appears as Question 27 in p.290 of
"Hyperidentities and Clones" \cite{DW}.

This presentation will be online with a URL like http://arxiv.org/abs/1408.XXXX . 

\newpage
\SetBgColor{red}
\SetBgOpacity{0.2}
\SetBgAngle{0}
\SetBgScale{50}
\SetBgContents{\textkorean{??}}
\SetBgHshift{5.5}
\BgThispage

\begin{center}
\section*{What Are Hyperidentities?}
\end{center}

Hyperidentities are statements in second-order
logic which can be viewed as an extension of
equational logic.  Informally, they are
universally quantified equations with the
 outer quantifiers over a domain of functions of certain
arities, and inner quantifiers over a domain of individuals.
The following example is written with subscripts to emphasize
that a symbol for a function variable such as $F_2$ and $G_3$
must have the same arity in all occurrences.
$$\forall F_2 \forall G_3 \forall xyz 
F_2(G_3(x,y,z),G_3(x,y,z)) \ideq G_3(F_2(x,x), F_2(y,y), F_2(z,z))$$

We will use a more casual style, suppressing quantifiers and subscripts.
We write the above as:
$$F(G(x,y,z),G(x,y,z)) \ideq G(F(x,x),F(y,y),F(z,z))$$

\textbf{What does this mean?}  Different things depending on the domain
of interpretation.  There will be one set $U$ of individuals and varying
sets of functions $\Omega$:

- Students of V. Belousov and Y. Movsisyan use the convention 
(e.g. in \cite{Be}, \cite{Mo} )
that $\Omega$ is a set of fundamental operations. This is in the
context of studying structures with finitely many operations $\cdot, \odot, \ldots$.
Some papers have these as binary quasigroup operations on the same set,
and the hyperidentity is a relation involving only these members of $\Omega$.

- Students of W. Neumann and W. Taylor (as in \cite{Ta}) will fix a similarity type $\tau$, say $\tau = <2,2,3,1>$ has symbols for two binary operations and a symbol
for a unary and another ternary operation.  $\Omega$ is then the collection of all
terms formed from these fundamental operations, including the projection operations
which appear as single variables.  This presentation will often use this
interpretation.  

- Students of K. Denecke, D. Schweigert, et.al. have a notion of
hypersubstitution which is covered in \cite{DW}; an effect is that certain
varieties of hyperidentities can be seen as having $\Omega$ range over certain
subsets of the set of terms.  In particular, prehyperidentities will exclude
the projection terms, and M-hyperidentities are a subset which are derived from a
monoid M of hypersubstitutions.  Parts of this presentation apply to this
interpretation.

\myfill

In all of these interpretations, arity is strictly preserved.  One informally
can substitute a projection function $p_i(x_1,\ldots,x_n) = x_i$ in 
interpreting an $n$-ary function symbol $F_n(x_1,...,x_n)$, but technically
one has to select only functions of the right arity $n$ from $\Omega$ for 
interpreting $F_n$.  

\myfill

A \textit{concrete clone} is collection of functions over an
algebra which contains the projection functions
and is closed under composition.  The hyperidentity
is interpreted as an equation between terms formed from 
clone members of certain arities.  Note that only Taylor's
use of hyperidentity corresponds to a clone identity.
The other uses correspond to the equation holding among
certain pairs of members of the clone, and not to all pairs
of given arities belonging to the clone.

\myfill 

\textbf{Technical Note:} To avoid complications in presentation,
we assume no constant symbols or functions of arity 0 , and we assume
the type has finitely many symbols of finite arity.  We call
such similarity types \textit{nice}.

Given a nice type $\tau$ and a hyperidentity $E$, we can form the
set of equations in the language of $\tau$-structures which come from 
all possible substitutions (sometimes restricted by terms which realize 
only functions in $\Omega$ for a given algebra) of a term of arity $n$ 
for the function symbol of same arity in $E$.  Given this set, it makes
sense to ask if it is logically equivalent (in equational logic) to
a finite set of equations in the same language.

\newpage
\SetBgColor{blue}
\SetBgContents{I}
\SetBgHshift{0}
\BgThispage

\begin{center}
\section*{A small basis for various hyperidentities}
\end{center}

We had been given the intuition that a small basis for hyperassociatvity
did exist.  We started by looking at consequences of small terms satisfying
the associative law.  We would use that to reduce the number of binary terms
to be considered.  $xy$ was a natural term to use, and let us
use semigroup notation and avoid many parentheses.
The starting term was $xyx$.

Ralph McKenzie pointed out to us that the semigroup variety given by 
the equation $xyxzxyx=xyzyx$ was locally finite.  
In particular, the 2-generated free semigroup of this variety was 
finite and we found it had 94 elements.  One could also use the
associative law for the term $xx$ ($x^2=x^4$) to reduce the number of terms 
to examine.

We then undertook to show that the remaining laws were a consequence
of two more laws $xxyyz=xxyxxyz$ and $xyyzz=xyzzyzz$.  Fortunately
this broke down into a few cases, based on the number of alternations
of $x$ and $y$ in the term $t$ ($x^ay^b, x^ay^bx^c, x^ay^bx^cy^d,$ and $xy^axyx$), the
last of which took advantage of analysis of smaller terms.  We did
the derivation, which was later used in Kunc's thesis \cite{Ku} .

\myfill
After getting this result, we did some more explorations and found
several other examples of pairs of hyperidentities and nice types which
gave finitely based equational theories.  

\myfill
Hyperidempotency: $$F(x) \ideq x,$$ finitely based for all types

\myfill
Hypermediality: $$F(G(x,y),G(z,w)) \ideq G(F(x,z),F(y,w)),$$
finitely based for all types containing function symbols of arity
at most 2,

\myfill
A version of entropic identity: 
$$F(G(x_11,\ldots,x_1n),\ldots,G(x_n1,\ldots,x_nn)) \ideq
G(F(x_11,\ldots,x_n1),\ldots,F(x_1n,\ldots,x_nn)),$$ finitely based
for types with symbols having arity at most $n$.

\myfill
HyperCommutativity: $$F(x,y)\ideq F(y,x),$$ finitely based for all types,
using the Taylor interpretation.  This included projections, yielding
the trivial base $x=y$.  For other interpretations which excluded
projections, work by Dudek and Kiesilewicz showed in \cite{KD}
that relative to the variety of semigroups, "totally commutative" 
semigroups belonged to one of four finitely based semigroup varieties.  
In unpublished work, McKenzie and Paseman showed that 
"total commutativity" was not finitely based relative to the variety 
of all magmas i.e. all universal algebras of type $<2>$.

\myfill
In the face of these examples, we conjectured that every hyperidentity
had a finitely based representation in a nice type.  We were wrong.
 
\newpage
\SetBgColor{red}
\SetBgContents{\textkorean{??}}
\SetBgHshift{5.75}
\BgThispage

\begin{center}
\section*{The Finite Basis Problem For Different Hyperidentities}
\end{center}

\myfill

Note that the finite basis problem is trivial when $\Omega$ is
finite.  The study of small bases can still be important, especially for
optimization problems.

\myfill

When $\Omega$ is the full term of clones, we have some degree of
uniformity.  We can use a model theory argument (sketched below) to show that
representing a finitely based hyperidentity is preserved across
reducts and a certain type of arity reduction.

\myfill

When one uses monoids of hypersubstitutions, one must take more care.
The sets of first order equations depends on the hypersubstitution,
and is usually restricted to within a particular type.  (There is also
some conflation between types used in the first order language and
the type of the language used to express the hyperidentity.  This is
necessary because the goal is to study things "in a first order way":
hypersubstitutions take place inside a monoid of terms of a given type,
and not as a substitution of algebraic operations for logical function
symbols.)  However, if the set of terms is "rich enough", one can
show an infinite basis in one type extends to an infinite basis in
a larger type using the same model theoretic argument.

\myfill

The basic argument involves infinite sequences of models.  Given
a sequence of $\tau$-models related to an equational theory derived from
a given hyperidentity $E$, we suppose this sequence witnesses the
fact that the theory is not finitely based, say by the $k$th model
satisfying only the first $k$ equations of the theory.  We create
a new sequence of $\sigma$-models, where $\sigma$ is like $\tau$, except
that $\sigma$ either has an extra function symbol, or $\sigma$ differs
from $\tau$ in that one function symbol has arity one more in $\sigma$
than in $\tau$.  It is then routine to construct a new sequence of $\sigma$-models
which witnesses the fact that the equational theory in the language $\sigma$
derived from $E$ is also not finitely based, by preserving as much of the
functionality of the given sequence of models.  We called this a "convexity
theorem"; the set of types which might not give rise to a finitely based
representation of a hyperidentity looked like a convex set in the partial order
of all nice types induced by the preserving relation.

\newpage
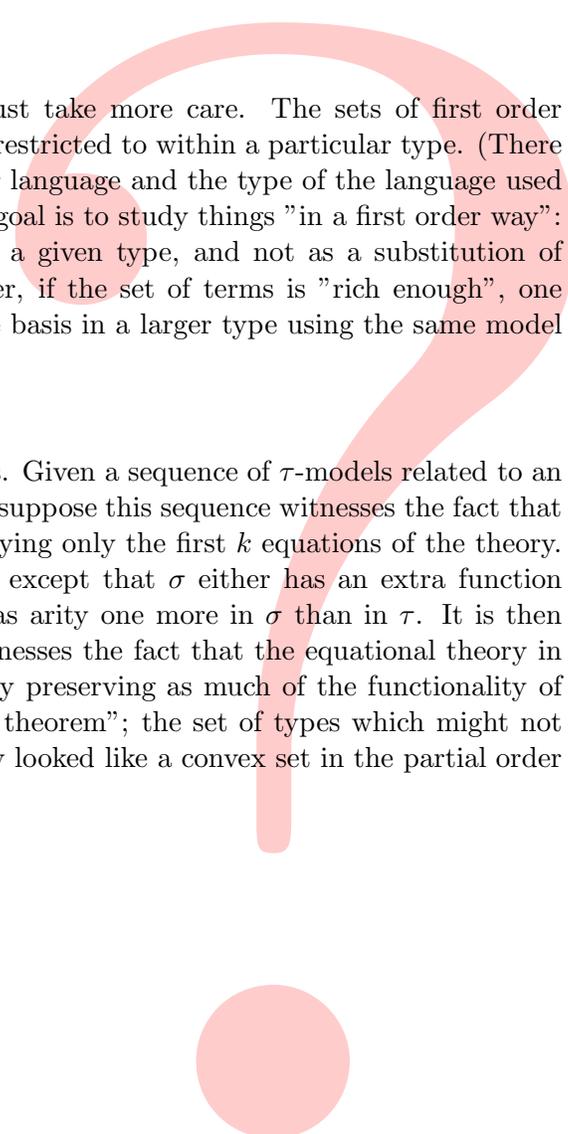
\SetBgColor{red}
\SetBgHshift{-6.25}
\BgThispage

\begin{center}
\section*{Going Back Twenty Years}
\end{center}

\myfill

We presented some of this work in \cite{P2} in Zurich ICM 1994.  Specifically:

\myfill

Hyperassociativity is finitely based for type $<2>$ (and many unary types)

\myfill

A list of various other $<$hyperidentity,type$>$ pairs yield finitely based
equational theories.

\myfill

$$F(F(x)) \ideq F(F(F(x)))$$ is not finitely based for the type $<1,1,1>$.
Our earlier conjecture that all such pairs were finitely based was wrong.

(The proof is based on the Thue-Morse word on three letters which is square
free.  We can choose certain long subwords of this and make large unary
terms $f(x)$ to substitute in the above hyperidentity.  The resulting term
$f(f(x))$ has no "small" subterms to use in a substitution, and so no "small"
equations can be used to derive $f(f(x))=f(f(f(x)))$ .  Thus there is no "small"
(finite) basis for the above theory.  We are indebted to Stuart Margolis for this 
suggestion.)

\myfill

There is a computable partial order $\prec$ on types which respects the property
of being finitely based:  if $\sigma \prec \tau$, then for any
hyperidentity $E$, $<E,\sigma>$ not finitely based implies $<E, \tau>$ is
also not finitely based.  (Conversely, $<E,\tau>$ is f.b. implies $<E,\sigma>$ is f.b. .
This comes from the "convexity theorem" mentioned in another panel.)

\myfill

The relation $\sigma \prec \tau$ includes and is generated by the relation $\sigma$ has
one less function symbol than $\tau$, and the relation $\sigma$ is the
same as $\tau$, except one symbol of $\tau$ has arity one greater than the
corresponding symbol of $\sigma$.

\myfill
Let there be a computable encoding of the countable set 
$$ A= \{ <E,\tau> \mid E \textrm{ is a hyperidentity and }
\tau \textrm{ is a nice type } \} . $$  
Let $B$ be that subset of $A$ for which there is
a finite equational basis which is logically equivalent to the first-order
equational theory.
The relation $\prec$ being computable implies that for fixed $E$, the set 
$B_E = \{<E,\tau> \in B \}$ is recursive.  Also, there is an algorithm to 
determine if the equivalent theory is trivial, i.e. has the basis $x=y$.

\myfill

Main Question: \textbf{Is B recursive?}

\newpage
\SetBgColor{blue}
\SetBgContents{C}
\SetBgHshift{0}
\BgThispage

\begin{center}
\section*{Hyperassociativity is Not Finitely Based for $<2,2>$}
\end{center}

\myfill 

This was one of the key results in \cite{P3}.  The preliminary
example showing $F^3(x) \ideq F^2(x)$ is not finitely based for
the type $<1,1,1>$ relied on finding an infinite sequence of unary
terms (in this case an infinite word), which would not allow for
appropriate substitutions to reduce large instances of the
hyperidentity to smaller instances.

Similarly, we needed an "infinite term" from which we could
use certain subterms to show that large instances of hyperassociativity
did not follow from smaller instances.  In this case the term was built up 
from alternating composition of one binary operation $\cdot$ with the
other $\circ$. The family of terms were used, along with a "dual" map
from terms which switched the binary operation symbols:

\begin{eqnarray*}
t_{00} & = & (x \cdot x) \\
t_{01} & = & (x \cdot y)  \\
t_{02} & = & (y \cdot x)  \\
t_{03} & = & (y \cdot y)    
\end{eqnarray*}
along with their *-counterparts, e.g. $t_{02}*=(y \circ x)$.  Then larger
terms were built using arrangements of the terms, as so:
\begin{eqnarray*}
t_{10} & = & t_{00}* \cdot  ( t_{01}* \cdot t_{02}* ) \\
t_{11} & = & t_{00}* \cdot  ( t_{01}* \cdot t_{03}* ) \\
t_{12} & = & t_{00}* \cdot  ( t_{02}* \cdot t_{03}* )  \\
t_{13} & = & t_{01}* \cdot  ( t_{02}* \cdot t_{03}* ) 
\end{eqnarray*}
along with their counterparts, e.g 
$t_{11}* = (x \cdot x) \circ ( (x \cdot y) \circ (y \cdot y))$.

The general pattern is given (for $n \geq 0$) by:
\begin{eqnarray*}
t_{(n+1)0} & = & t_{n0}* \cdot  ( t_{n1}* \cdot t_{n2}* ) \\
t_{(n+1)1} & = & t_{n0}* \cdot  ( t_{n1}* \cdot t_{n3}* ) \\
t_{(n+1)2} & = & t_{n0}* \cdot  ( t_{n2}* \cdot t_{n3}* ) \\
t_{(n+1)3} & = & t_{n1}* \cdot  ( t_{n2}* \cdot t_{n3}* ) \\
t_{(n+1)0}* & = & t_{n0} \circ  ( t_{n1} \circ t_{n2} )   \\
t_{(n+1)1}* & = & t_{n0} \circ  ( t_{n1} \circ t_{n3} )   \\
t_{(n+1)2}* & = & t_{n0} \circ  ( t_{n2} \circ t_{n3} )   \\
t_{(n+1)3}* & = & t_{n1} \circ  ( t_{n2} \circ t_{n3} ) 
\end{eqnarray*}
Then  $t=t_{(2k+1)0}(x,t_{(2k+1)0}(y,z))$
has no "small" substitution instances
appearing inside it.  One
would need such instances if there were a finite basis for
hyperassociativity to transform that term to
$u=t_{(2k+1)0}(t_{(2k+1)0}(x,y),z))$, but the only such allowed terms
have more than $k$ binary operation symbols, which means one needs
large instances of hyperassociativity to prove this instance $t=u$
in equational logic.

\myfill

Similar arguments were used to show hyperassociativity is not finitely
based for the types $<3>$ and for $<2,1>$.  This establishes the finite
basis character of hyperassociativity for all nice types.

\newpage
\SetBgColor{red}
\SetBgContents{\textkorean{??}}
\SetBgHshift{5.5}
\BgThispage
\begin{center}
\section*{Aside: Hypersubstitutions and model theory}
\end{center}

\myfill

This portion diverges from considering the two problems,
but is of related interest.

\myfill

\cite{DK} provided the inspiration for the original problem
on hyperassociativity and some results in \cite{P1}.  However,
the desire to understand the model theoretic side spurred
the development arc in \cite{P1}, \cite{P2}, \cite{P3}, and
\cite{P4}.

\myfill
Corresponding to the linguistic side of hyperidentities is
the model theoretic side of hypervarieties and solid varieties. 
Hypervarieties according to \cite{Ta} are "varieties of varieties",
or classes of varieties which are closed under certain operations
on varieties such as varietal product and reduct.  While of some
interest, the study of clone identities for abstract and concrete
clones has appeared more than studies focussed on hypervarieties.

\myfill
Solid varieties are the by-product of the research direction using
hypersubstitutions.  Briefly, one considers maps on the term algebra
of a type $\tau$ to itself, and one can form a monoid of such maps
which stand for replacing a basic function symbol $\tau$ of arity 
$n$ by a term from $\tau$ of the same arity.  One then sees a
hyperidentity as relation among hypersubstitutions, and one can
use this idea to form derived algebras, which are algebras enhanced
with term operations arising from these hypersubstitutions.  A solid
variety is now a variety containing these derived algebras, and the
collection of solid varieties of  $\tau$-algebras forms
a complete sub-lattice of the lattice of varieties of all $\tau$-algebras.
As Polak showed in \cite{P9} even for semigroups there is some richness
to this sublattice.

\myfill

The Finite Basis Problem should relate to locally finite
solid varieties, and hopefully reveal their relation
to other solid varieties.  In particular it should help
talk about doing computations (or not) in the lattice
of solid varieties,
perhaps determine whether such a variety represents a 
join-irreducible element, or reveal other information.
\myfill

The Finite Basis Problem should also play a role in
seeing how complex the lattice can be in the case 
the hyperidentity is not finitely based with respect to
a given type.  Further, the preorder $\prec$ should give
rise to a categorical mapping between solid varieties
of $\prec$-related types, possibly providing embeddings
of one lattice of solid varieties into another.

\newpage
\SetBgHshift{-6.25}
\BgThispage
\begin{center}
\section*{A Theory of Small Trees, With Applications}
\end{center}

\myfill

In \cite{P4}, we reported some of the progress in finding
non-finitely based pairs, particularly those involving
hyperassociativity and the types  $<2,2>$,
$<2,1>$, and $<3>$.

\myfill

Also, Ralph McKenzie and George Bergman had shown us similar
families of terms witnessing that hypermediality was not finitely
based for some types.

\myfill

Finally, as the first example  of a nfb representation came
from symbolic dynamics, we suggested in \cite{P4}
the notion that symbolic dynamics should be
extended to solve the Finite Basis problem.  

\myfill

In particular, instead of a tree representation for a term,
consider a class of decorated directed sets, or infinite trees.
(every two nodes belong to some common infinite sub term, each
node decorated with a function symbol whose arity corresponded
to the number of children of that node.).  Such an infinite
tree could be used to derive or otherwise represent
an infinite sequence of terms, 
with earlier terms being subterms of later terms. One now has
the concept of which trees map to their own shifts, much as
occurs in symbolic dynamics.  Indeed, standard symbolic dynamics
on $\tau^{\textrm{Z}}$ (where $\tau$ is playing a role as both
a type of unary functions and a symbol alphabet) is a special
case of this, the infinite term of composed unary functions
corresponding to a trajectory.

\myfill

One can study collections of such trees to see how they act
under generalized shifts.  In particular, where certain
subterms correspond to one side of a hyperidentity, one might
find rigid trees which have few or no such shifts that would
correspond to substitutions in a derivation in equational
logic.  Such a "rigid" tree would serve a role analogous
to the infinite sequence of models in witnessing the lack
of a possibility for a finite basis for the hyperidentity.

\myfill
In addition to finding collections of "rigid" trees, one
needs to determine which aspects of these collections lend
to computability results.  If it turns out that 
"rigid" trees form a non-recursive set, that may lead to
a negative answer of the Finite Basis Problem.

\myfill

At the time of this presentation, we were becoming aware of the
work following Denecke's school, but did not see how the concept
of hypersubstitution was changing the face of the research on
hyperidentities.  We think that some computability issues raised
in \cite{P4} could inspire algorithmic complexity issues for
hyperidentities.

\newpage
\SetBgColor{blue}
\SetBgContents{M}
\SetBgHshift{0}
\BgThispage
\begin{center}
\section*{The case $xxyyz = xxyxxyz$}
\end{center}

\myfill

If we consider a variation of the argument for the hyperidentity
$F(F(x)) \ideq F(F(F(x)))$, we should find in the semigroup variety
satisfying $xxyyz=xxyxxyz$ relatively-free 3 or 4-generated free
algebras which are not finite.  This is because we should find
semigroup words in 3 or 4 variables which are squarefree, and will
be distinct in the free algebra of this variety.

\myfill

Even though we do not have the local finiteness of semigroups in the case
$xyxzxyx=xyzyx$, it is enough to have the relatively free 2-generated
semigroup be finite.  This amounts to looking for repeated squares
in words on an alphabet of two letters.  Fortunately there are enough
repeated squares that we can claim a similar finiteness result 

\myfill

Note that we have $x^2y^2x^2z=x^2y^2x^4z$, and $x^2y^bz=(x^2y)^bz$.
Further, let $W(x,y)$ and $Z(x,y)$ be words in $x$ and $y$, with $Z(x,y)$
not the empty word.  
By the two consequences above, any word in $x$ and $y$ which is 
 of the form $W(x,y)x^2y^3Z(x,y)$ or $W(x,y)x^2y^2x^2Z(x,y)$ has a logical
equivalent where $Z(x,y)$ is replaced by $Z'(x,y)$.  $Z'(x,y)$
has all but the last power reduced modulo 2, e.g. if
$Z(x,y)$ is $x^{a(1)}y^{a(2)}\ldots x^{a(2n-1)}y^{a(2n)}$, and
$b(i) = a(i) \bmod 2$, and one writes the empty string as $x^0$ or $y^0$,
then $Z'(x,y)$ is $x^{b(1)}y^{b(2)}\ldots x^{b(2n-1)}y^{a(2n)}$.
This combined with $x^2y^3z=x^2y^5z$ gives finitely many alternatives
for $Z'(x,z)$.

\myfill

Similarly, one can reduce large patterns strictly containing $xyxyZ(x,y)^2$.
One then finds only finitely many inequivalent alternatives for $W(x,y)$.
This results in a finite set of words, allowing another route to showing
that hyperassociativity is finitely based for $<2>$.  

\myfill

Of course, one can perform a
similar analysis with the subvariety which also satisfies $xyyzz=xyzzyzz$.
This addresses problem 28 from \cite{DW}.
\newpage
\SetBgColor{green}
\SetBgContents{2}
\SetBgHshift{2}
\BgThispage
\begin{center}
\section*{Future Directions and Applications}
\end{center}

\myfill

In unpublished work, we found that an analogue to Burnside's problem
was not finitely based.  This means that $F^n(x)\ideq F(x)$ was not
finitely based for the type $<2>$, and so for "larger" types.

\myfill
One unexplored avenue regards categorical interpretations of the
problem.  The clones correspond to algebraic theories of Lawvere,
which have a substantial literature.  We are unsure how the 
Finite Basis problem is realized in this context.

\myfill
At the moment, the connection to symbolic dynamics seems intriguing.
We hope to extend the notion of dynamics on $S^{Z}$ to dynamics on
directed sets, sections of which represent "infinite terms".  We hope
this extension of dynamics will allow the analysis and solution of the Main Question.

\myfill
We also note work of Melkonian \cite{Me} in circuit optimization, who used 
hyperidentities and hypergates, to analyze and solve a k-out-of-n circuits
problem.  We hope to see more applications to circuit optimization.

\newpage
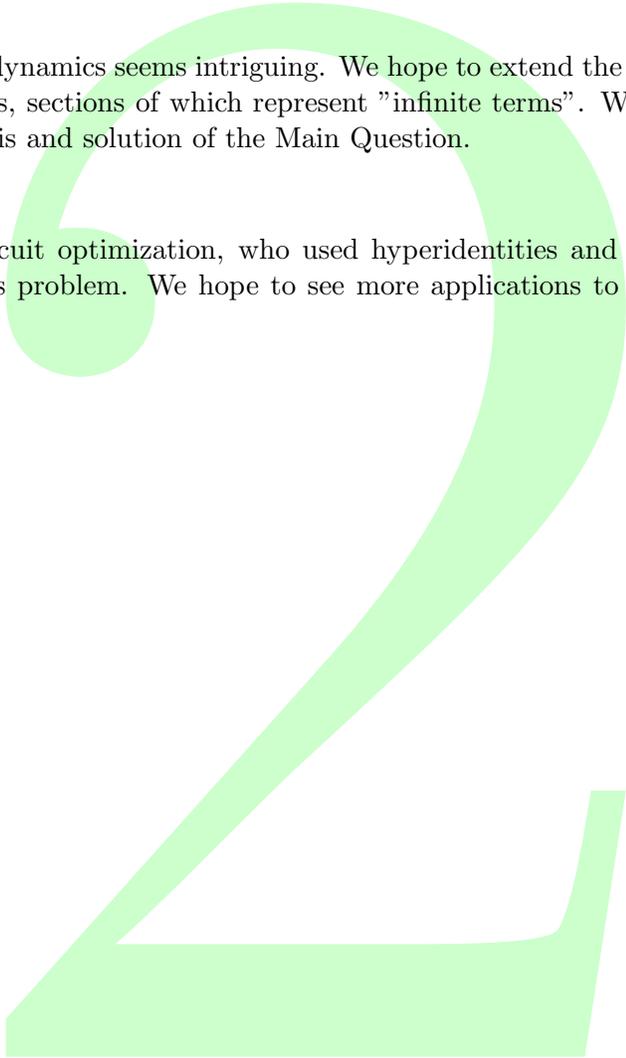
\SetBgContents{01}
\SetBgHshift{0}
\BgThispage
\begin{center}
\section*{Acknowledgments}
\end{center}
\myfill

I wish to acknowledge my advisors at U.C. Berkeley, Ralph McKenzie and
George Bergman.  

\myfill
I also appreciate others who encouraged me, especially Stuart Margolis, and
 including
the attendees and speakers of the General Algebra Seminar, among them
Art Drisko, Jonathan Farley, Rad Dimitric, Charles Latting, Andrew Ensor,
Japheth Wood, Arturo Magidin, Stefan Schmidt, J.D. Phillips, and William Craig.

\myfill
I appreciate the people who gave me opportunities to present and study
this subject and who cited my unpublished work.  These include
  Dietmar Schweigert, Klaus Denecke, Shelly Wismath,
Libor Polak, and Yuri Movsisyan.  Also their students and colleagues
among them Jorg Koppitz and Michal Kunc.

\myfill
I thank the Organizing Committee for the Congress for giving me the
opportunity to present this work, and the audience for reading it.

\myfill
I thank Joy Song for her help on some aspects of the poster design.

\myfill
I am grateful for the love and support of my family, including my wife Cheryl
and son Richard.

\newpage
\SetBgContents{4}
\SetBgHshift{-2}
\BgThispage

\end{document}